\documentclass{article}

\usepackage{amssymb}
\newtheorem{satz}{Satz}[section]
\newtheorem{remark}[satz]{Remark}
\def\proofbegin{{\textit  Proof}. }
\newtheorem{thm}{Theorem}
\newtheorem{lem}[thm]{Lemma}
\newtheorem{defn}{Definition}
\title{Numerical Computation of Takens-Bogdanov Points for Delay Differential Equations\thanks{Supported by NSFC grants 10971022 and 11071102.}}

\author{{ Yingxiang Xu}\footnote{E-mail: yxxu@nenu.edu.cn} and { Vital D.~Mabonzo}\\[3mm] {School of Mathematics and Statistics}, \\{Northeast Normal University, Changchun 130024, China.}}

\date{}

\begin{document}
\maketitle
\begin{abstract}
The paper presents a numerical technique for computing directly the
Takens-Bogdanov points in the nonlinear system of differential
equations with one constant delay and two parameters. By representing
the delay differential equations as abstract ordinary differential
equations in their phase spaces, the quadratic Takens-Bogdanov point is
defined and a defining system for it is produced. Based on the descriptions for the eigenspace associated with the double zero eigenvalue, we reduce the defining system to a finite dimensional algebraic equation. The quadratic Takens-Bogdanov point, together with the corresponding values of parameters, is proved to be the regular solution of the reduced defining system and then can be approximated by
the standard Newton iteration directly.

\end{abstract}


\pagestyle{myheadings} \markboth{Y. Xu et al. Applied Mathematics and Computation xxx (2010) xxx-xxx}{Numerical computation of T-B points for DDEs}

\section{Introduction}
In this paper, we discuss the computation of the Takens-Bogdanov points in the following delay differential equations (DDEs)
\begin{equation}\label{eintro1}
\dot x(t)=f(x(t),x(t-\tau),\lambda,\mu),
\end{equation}
where $\lambda$, $\mu\in\mathbb{R}$ are parameters, $\tau>0$ is a constant delay. The DDEs have been widely studied  because they often give more accurate descriptions for the phenomena in nature and engineering by taking into account not only the present state but also their histories (see \cite{bc63}, \cite{dr77}, \cite{hl93}, \cite{ahd06} etc.).

Takens-Bogdanov (short for T-B) bifurcation is one of the important bifurcations in dynamical systems, it explains the mechanics for the occurrence of the Hopf point branch and the homoclinic branch as well as the saddle-node point branch. It acts as a bridge for investigating the global dynamical behavior through the local properties of the dynamical systems. Theoretically speaking, based on the local bifurcations near T-B points, one can obtain the T-B points on the paths of fold points (when a second real eigenvalue crosses the imaginary axis) or on the paths of Hopf points (when the two conjugate pure imaginary eigenvalues coalesce). For details, one can refer to \cite{ts74}, \cite{bv75} for ODEs and to \cite{fa95a}, \cite{xh08} for DDEs.

The paths of fold points or Hopf points are generally traced numerically by the continuation method \cite{Ag94}. However, when T-B point is encountered along the paths, the continuation technique is not valid any more. As a result, to compute directly the T-B point is of great importance. In addition, it is the foundation of branch switching, that is to compute the other solution branch, i.e., homoclinic branch or another solution branch differing from the paths traced before, emanating from T-B point. In general, the T-B point is approximated by numerical iteration method applied to an enlarged system. In this process, an initial value is necessary and can be obtained by detecting the changes of eigenvalues along the paths of fold points or Hopf points by the continuation technique \cite{g00}. For DDEs, the work could be done by the method developed by Luzyanina and Roose \cite{Dr96}, where a defining system is given to determine the Hopf bifurcation point algebraically, and the continuation techniques are considered. These techniques are enclosed in DDE-BIFTOOL \cite{elr02}. Of course, to start the continuation of the path of fixed points, a good approximation is needed. In fact, in terms of that the constant delays do not affect the position of equilibria of the system, the continuation methods for the equilibria of ODEs could be employed to DDEs \cite{Dr96, rs07} directly.  Unfortunately, noting that T-B point is a singularity of codimension 2, the techniques for computing the Hopf points developed in \cite{Dr96} for DDEs could not be applied to compute the T-B points directly. Besides, many applications require us to obtain the T-B points numerically first so that we can start our continuation from the T-B point (e.g. \cite{fgks}). In addition, tracing the T-B points branch emanating from the singular point of higher order needs to compute T-B points first as well. Therefore, to determine the T-B points numerically is an important task in the numerical analysis for dynamical systems.

It is well studied for ODEs to calculate the T-B points, one can refer to Griewank and Reddien \cite{Gr89}, Govaerts
\cite{gw93}, Beyn \cite{Beyn94}, Yang \cite{zh07} and the subsequent
articles by many authors. The method is described sketchily as follows.

 Consider the following parameterized dynamical systems
\begin{equation}\label{eyan}
\dot x=f(x,\lambda,\mu),\quad x\in\mathbb{R}^n,\quad
\end{equation}
where $\lambda,\mu\in\mathbb{R}$ are parameters,
$f:\mathbb{R}^n\times\mathbb{R}\times\mathbb{R}\rightarrow\mathbb{R}^n$
is a continuously differentiable function with
$f(x^0,\lambda^0,\mu^0)=0$, i.e. $ x^0$ is an equilibrium of
(\ref{eyan}) for $\lambda=\lambda^0$ and $\mu=\mu^0$.

\begin{defn}{\rm\cite{zh07}}\label{def1}
$(x^0,\lambda^0,\mu^0)$ is called a T-B point of (\ref{eyan}) if
the following conditions hold:
\begin{description} \item{\rm(i)}  $f(x^0,\lambda^0,\mu^0)=0$;
\item{\rm(ii)}  $\mathcal{N}(f_x^0)={\rm span}\{\eta_0\}$,
$\eta_0\neq 0$; \item{\rm(iii)}
$\mathcal{R}(f_x^0)=\{\varsigma\in\mathbb{R}^n,\xi_0^T\varsigma=0\}$;
\item {\rm(iv)} $\eta_0\in \mathcal{R}(f_x^0)$.
\end{description}
\end{defn}
The conditions in Definition \ref{def1} imply that there exist $\eta_1\in M$ and $\xi_1\in
M_1$, such that
\[
\left\{\begin{array}{l}f_x^0\eta_1+\eta_0=0,\\
l_0^T\eta_1=0,\end{array}\right.
\]
and
\[
\left\{\begin{array}{l}(f_x^0)^T\xi_1+\xi_0=0,\\
\xi_1^T\omega_0=0,\end{array}\right.
\]
where $\mathbb{R}^n=\mathcal{N}(f_x^0)\oplus M$,
$\mathbb{R}^n=\mathcal{N}((f_x^0)^T)\oplus M_1$, $l_0^T\eta_1=0$
is equivalent to $\eta_1\in M$, $\xi_1^T\omega_0=0$ is equivalent
to $\xi_1\in M_1$.

\begin{defn}{\rm\cite{zh07}}\label{d3}
$(x^0,\lambda^0,\mu^0)$ is called a quadratic T-B point of
(\ref{eyan}), if it is a T-B point and\\
{\rm(i)} $\xi_0^Tf_\lambda^0\neq 0;$\\
{\rm(ii)}
$\bar d_0=\det\left(\begin{array}{ll}\xi_0^T\bar A\eta_0&\xi_0^T\bar B\eta_0\\[3mm]
\begin{array}{l}\xi_0^T\bar A\eta_1+\xi_1^T\bar A\eta_0
\end{array}
&
\begin{array}{l}\xi_0^T\bar B\eta_1+\xi_1^T\bar B\eta_0\end{array}\end{array}\right)\neq
0$;\\
{\rm(iii)} $\xi_0^T\eta_1\neq 0,$\\
where $\bar A=f_{xx}^0\eta_0$, $\bar
B=f_{xx}^0\bar\nu_{\lambda\mu}+\bar
c_{\lambda\mu}f_{x\lambda}^0+f_{x\mu}^0$ with $\bar
c_{\lambda\mu}=-\frac{\xi_0^Tf_\mu^0}{\xi_0^Tf_\lambda^0}$ and
$\bar \nu_{\lambda\mu}$ satisfying
\[\begin{array}{l}
f_x^0\bar\nu_{\lambda\mu}+\bar c_{\lambda\mu}f_\lambda^0+f_{\mu}^0=0,\\
l_0^T\bar \nu_{\lambda\mu}=0.
\end{array}\]
\end{defn}
The defining system for the quadratic T-B points of (\ref{eyan}) is then
given by
\begin{equation}\label{ed7}
H_1(w_1)=\left(\begin{array}{l}f(x,\lambda,\mu)\\
f_x(x,\lambda,\mu)\eta\\
f_x(x,\lambda,\mu)\zeta+\eta\\
l_0^T\eta-1\\
l_0^T\zeta\end{array}\right)=0,
\end{equation}
where $w_1=(x,\eta,\zeta,\lambda,\mu)^T\in
U=\mathbb{R}^n\times\mathbb{R}^n\times\mathbb{R}^n\times\mathbb{R}\times\mathbb{R}$,
$H_1:U\rightarrow
U$.

\begin{thm}{\rm\cite{zh07}}\label{tf}
Assume $(x^0,\lambda^0,\mu^0)$ is a quadratic T-B point of
(\ref{eyan}), then $H_1(w_1)$ is regular at $w_1^0=(x^0,\eta_0,\eta_1,\lambda^0,\mu^0)^T$.
\end{thm}

By Theorem \ref{tf}, the quadratic T-B point of (\ref{eyan}) can be obtained
numerically by applying iteration methods to the defining system
(\ref{ed7}). However, unlike (\ref{ed7}), the defining system for T-B point of (\ref{eintro1}) should be defined in a Banach space from the viewpoint of \cite{Gr89}, see section 2 for details.  If one wants to solve such a system directly, certain discretization must be applied and the discretization error cannot be avoided.

In this paper, we present a numerical technique for computing the T-B points in nonlinear system of DDEs with one constant delay and two parameters, which carries the results of computing the T-B points of ODEs into the case of DDEs.  Based on the descriptions for the eigenspace associated with the double zero eigenvalue, we reduce the defining system to a nonlinear algebraic equation with finite dimension to avoid discretizing it directly.  The quadratic T-B point, together with certain values of the parameters, is proved to be the regular solution of the reduced defining system. Therefore, the quadratic T-B point of DDEs can be solved numerically by the classical iteration methods.

The paper is arranged as follows. Through representing the DDE as an abstract ODE in Section 2, we define the quadratic T-B point and produce a defining system for the quadratic T-B point of DDEs according to the methods described above. In Section 3, by using of the descriptions for the eigenspace associated with double zero eigenvalue, we simplify the definition of the quadratic T-B point of DDEs, meanwhile, the defining system obtained in Section 2 is reduced to a finite dimensional one, which is proved to be regular at the quadratic T-B points such that can be solved by the standard Newton iteration. To show the efficiency of our method, a numerical experiment is carried for a predator prey system with time delay in the last section.

\section{Regular defining system for the T-B points in Banach space}
By a simple change of timescale, one may take $\tau=1$ in (\ref{eintro1}). Hence we consider the following DDEs
\begin{equation}\label{eyan1}
\dot x(t)=f(x(t),x(t-1),\lambda,\mu),
\end{equation}
where $\lambda$,$\mu\in\mathbb{R}$ are parameters,
$x\in\mathbb{R}^n$, $f(x,y,\lambda,\mu)$ is a $C^r(r\geq2)$ smooth
function from
$\mathbb{R}^n\times\mathbb{R}^n\times\mathbb{R}\times\mathbb{R}$\quad
to $\mathbb{R}^n$. Assume $x^0$ is the equilibrium of
$f(x(t),x(t-1),\lambda,\mu)$ at parameter value
$\lambda=\lambda^0$ and $\mu=\mu^0$. Eq. (\ref{eyan1}) can be
linearized at $(x^0,\lambda^0,\mu^0)$ as
\begin{equation}\label{elinear}
\dot x(t)=f_1^0(x(t)-x^0)+f_2^0(x(t-1)-x^0),
\end{equation}
where $f_1^0=\frac{\partial f}{\partial
x}(x^0,x^0,\lambda^0,\mu^0)$, $f_2^0=\frac{\partial f}{\partial
y}(x^0,x^0,\lambda^0,\mu^0)$.

Denote by $C:=C([-1,0],\mathbb{R}^n)$ the Banach space  of continuous
mappings from $[-1,0]$ to $\mathbb{R}^n$ with norm
$\|\phi\|=\max_{\theta\in[-1,0]}|\phi(\theta)|$ ($|\cdot|$ is a some
norm in $\mathbb{R}^n$), and let
$$
\eta_{\lambda,\mu}(\theta) = \left\{ \begin{array}{ll}
f_1(x^0,x^0,\lambda,\mu)+f_2(x^0,x^0,\lambda,\mu),&\theta =0 ,\\
f_2(x^0,x^0,\lambda,\mu),&\theta \in (-1,0),\\
0,&\theta =-1
\end{array}
\right.
$$ be a bounded variation matrix-valued function on $[-1,0]$ parameterized by $\lambda$ and $\mu$.
Noting that
$$\int_{-1}^0d\eta_{\lambda,\mu}(\theta)(x_t(\theta)-x^0)=f_1(x^0,x^0,\lambda,\mu)(x(t)-x^0)+f_2(x^0,x^0,\lambda,\mu)(x(t-1)-x^0),$$
where $x_t(\theta)=x(t+\theta)$, the linear operator from $C$ to
$\mathbb{R}^n$ defined by
$$L_{\lambda,\mu}x_t=\int_{-1}^0d\eta_{\lambda,\mu}(\theta)x_t(\theta)$$ is bounded.

The solutions of (\ref{elinear}) generate a $C_0$-semigroup
$\{T_0(t),t\geqslant 0\}$ on $C$ with infinitesimal generator
$\mathcal{A}_0:C\longrightarrow C$ defined by
\[\mathcal{A}_0\phi=\dot\phi,\]
\[
D(\mathcal{A}_0)=\begin{matrix}\{\phi\in C^1([-1,0],\mathbb{R}^n);
\dot\phi(0)=\int_{-1}^0d\eta_{\lambda^0,\mu^0}(\theta)\phi(\theta)\}\end{matrix}.
\]
It is known that the spectrum of the operator $\mathcal{A}_0$ consists
of its point spectrum only \cite{hl93}, i.e.
$\sigma(\mathcal{A}_0)=\sigma_p(\mathcal{A}_0)$, and
$z\in\sigma_p(\mathcal{A}_0)$ if and only if
\begin{equation}\label{character}
\det(z I-f_1^0-f_2^0e^{-z})=0.
\end{equation}
Usually we call (\ref{character}) the characteristic equation of
(\ref{elinear}).

To develop  the theory for computing the T-B points of DDEs, it is convenient to write
(\ref{eyan1}) as the following abstract ODE in $C$ \cite{hl93}
\begin{eqnarray}\label{eyan5}
\frac{d}{dt}u=G(u,\lambda,\mu),
\end{eqnarray}
where
\begin{eqnarray}\label{eyan6}
G(u,\lambda,\mu)(\theta) = \left\{ \begin{array}{ll}
f(u(0),u(-1),\lambda,\mu),&\theta =0, \\
\dot u(\theta),&\theta \in [-1,0]
\end{array}
\right.
\end{eqnarray}
and the domain of $G(\cdot,\lambda,\mu)$ is $\{u\in C^1:\ \dot
u(0)=f(u(0),u(-1),\lambda,\mu)\}.$

 Let $C^\ast=C([0,1],\mathbb{R}^{n\ast})$ be the adjoint space of
 $C$, with $\mathbb{R}^{n\ast}$ the $n$-dimensional space of row vectors. The adjoint bilinear form on
 $C^\ast\times C$ is defined by \cite{hl93}
 $$
(\psi,\phi)=\psi(0)\phi(0)-\int_{-1}^0\int_{0}^\theta\psi(\xi-\theta)d\eta_{\lambda^0,\mu^0}(\theta)\phi(\xi)d\xi.
 $$
The T-B point of (\ref{eyan1}) can be determined through investigating the T-B point of (\ref{eyan5}) since they share the same equilibria. Here we define the T-B point of (\ref{eyan5}) first.
 \begin{defn}\label{dtdde}
$(u^0,\lambda^0,\mu^0)$ is called a T-B point of (\ref{eyan5}), or equivalently of (\ref{eyan1}), if\\
{\rm(i)} $G(u^0,\lambda^0,\mu^0)=0;$\\
{\rm(ii)} $\mathcal{N}(G_u^0)={\rm span}\{\phi_1\}, \phi_1\neq 0$;\\
{\rm(iii)} $\mathcal{R}(G_u^0)=\{v\in C, (\psi_2, v)=0\}$;\\
{\rm(iv)} $\phi_1\in \mathcal{R}(G_u^0)$.
\end{defn}
The condition (iv) in Definition \ref{dtdde} is equivalent to
$(\psi_2,\phi_1)=0$, and there exist $\psi_1\in C^\ast$ and
$\phi_2\in C$ such that
\[\begin{array}{ll}
(G_u^0)^\ast\psi_1+\psi_2=0,& \ \ (\psi_1,\tilde\omega_0)=0,\\
G_u^0\phi_2+\phi_1=0,&\ \ (\tilde l_0,\phi_2)=0, \end{array}\]
with $\tilde\omega_0\in C$ and $\tilde l_0\in C^\ast$ satisfying
\[
(\phi_2,\tilde\omega_0)-1=0,\ \ (\tilde l_0,\phi_1)-1=0.
\]
\begin{defn}\label{d2.2}
$(u^0,\lambda^0,\mu^0)$ is called a quadratic T-B point of
(\ref{eyan5}), or equivalently of (\ref{eyan1}), if it is a
T-B point and\\
{\rm(i)} $(\psi_2,G_\lambda^0)\neq 0$;\\
{\rm(ii)} $d_0=\det \left(\begin{array}{cc} (\psi_2,
A\phi_1)&(\psi_2,
B\phi_1)\\
(\psi_2,A\phi_2)+(\psi_1,A\phi_1)&(\psi_2, B \phi_2)+(\psi_1,
B\phi_1)\end{array}\right)\neq 0;$\\
{\rm(iii)}$(\psi_2,\phi_2)\neq 0,$\\
where $A=G_{uu}^0\phi_1$,
$B=G_{uu}^0\nu_{\lambda\mu}+c_{\lambda\mu}G_{u\lambda}^0+G_{u\mu}^0$
with $c_{\lambda\mu}=-\frac{(\psi_2,G_\mu^0)}{(\psi_2,G_\lambda^0)}$
and $\nu_{\lambda\mu}$ satisfying
\[\begin{array}{l}
G_u^0\nu_{\lambda\mu}+c_{\lambda\mu}G_\lambda^0+G_{\mu}^0=0,\\
(\tilde l_0,\nu_{\lambda\mu})=0.\end{array}\]
\end{defn}
Based on the preparations above, we produce the following defining
system for the quadratic T-B points of (\ref{eyan1})
\begin{equation}\label{ea2}
H_2(w_2)=\left(\begin{array}{l} G(u,\lambda,\mu)\\
G_u(u,\lambda,\mu)e_1\\
G_u(u,\lambda,\mu)e_2+e_1\\
(\tilde l_0,e_1)-1\\
(\tilde l_0,e_2)
\end{array}\right)=0,
 \end{equation}
where $w_2=(u,e_1,e_2,\lambda,\mu)^T\in Y=C\times C\times C\times
\mathbb{R}\times\mathbb{R}$, $H_2:Y\rightarrow
Y$.
\begin{thm}Assume $(u^0,\lambda^0,\mu^0)$ is a quadratic T-B point of
(\ref{eyan5}), or equivalently of (\ref{eyan1}), then $H_2(w_2)$ is
regular at
$w_2^0=(u^0,\phi_1,\phi_2,\lambda^0,\mu^0)^T$.
\end{thm}

\proofbegin
The proof is quite similar to Theorem \ref{tf}.
$\Box$
\section{Simplification  of the defining system}
The defining system (\ref{ea2}) is defined in a Banach space, so if we solve it directly, we will encounter many difficulties. The space $C$ must be discretized first whenever what kinds of numerical methods are applied. This needs to store a large amount of data and causes the discretization error undoubtedly, which is not expected surely. Another big difficulty lies in the form of
the function $G(u,\lambda,\mu)$ which needs to be dealt with piecewisely
since its special form caused by representing the DDE as an abstract ODE.

In this section, we will reduce (\ref{ea2}) to an equivalent form, which has finite
dimension and can be solved easily.

We first give an equivalent definition for the T-B
points of (\ref{eyan1}) by the infinitesimal generator
$\mathcal{A}_0$ of the $C_0$-semigroup defined by the solutions of
(\ref{elinear}).
\begin{defn}\label{dtdde2}
$(x^0,\lambda^0,\mu^0)$ is called a T-B point of (\ref{eyan1}), if\\
{\rm(i)} $f(x^0,x^0,\lambda^0,\mu^0)=0;$\\
{\rm(ii)} $\mathcal{N}(\mathcal{A}_0)={\rm span}\{\phi_1\}, \phi_1\neq 0$;\\
{\rm(iii)} $\mathcal{R}(\mathcal{A}_0)=\{\gamma\in C, (\psi_2, \gamma)=0\}$;\\
{\rm(iv)} $\phi_1\in \mathcal{R}(\mathcal{A}_0)$.
\end{defn}
The condition (iv) in Definition \ref{dtdde2} is equivalent to
$(\psi_2,\phi_1)=0$, and there exist $\phi_2\in \bar M$ and
$\psi_1\in \bar M_1$ such that
\[\begin{array}{ll}
\mathcal{A}_0\phi_2+\phi_1=0,& \ (l,\phi_2)=0,\\
\mathcal{A}_0^\ast\psi_1+\psi_2=0, &\ (\psi_1,\varpi)=0,
\end{array}\]
where $C=\mathcal{N}(\mathcal{A}_0)\oplus\bar M$,
$C=\mathcal{N}(\mathcal{A}_0^\ast)\oplus\bar M_1$, $(l,\phi_2)=0$
is equivalent to $\phi_2\in\bar M$ with $l=l_1-sl_2$, $s\in[0,1]$;
$(\psi_1,\varpi)=0$ is equivalent to $\psi_1\in\bar M_1$.

Noting that $x^0\in\mathbb{R}^n$ can be regarded as a special
element in $C$ coincided with $u^0$,
$(x^0,\lambda^0,\mu^0)$ is a quadratic T-B point if it satisfies
Definition \ref{dtdde2} and \ref{d2.2}. The following theorem can be used to determine whether $(x^0,\lambda^0,\mu^0)$ is a T-B point of (\ref{eyan1}) or not.

\begin{thm}{\rm\cite{xh08}}\label{t1}
Assume ${\rm Re}{z}\neq0$ if $z\in
\sigma_p(\mathcal{A}_0)\backslash \{0\}$, Eq. (\ref{eyan1}) has a
T-B singularity if and only if the following conditions hold:
\[\label{e2.6}
\begin{array}{l}
(i)\ \ {\rm rank}( f_1^0+f_2^0)=n-1;\\
(ii)\ \ \mbox{if }\mathcal{N}(f_1^0+f_2^0)={\rm span}\{\phi_1^0\},
\mbox{ then }
 (f_2^0+I)\phi_1^0\in {\mathcal{R}}(f_1^0+f_2^0);\\
(iii)\ \ \mbox{if }(f_1^0+f_2^0)\phi_2^0=(f_2^0+I)\phi_1^0,\mbox{
then } (f_2^0+I)\phi_2^0-\frac{1}{2}f_2^0\phi_1^0\not\in
 {\mathcal{R}}(f_1^0+f_2^0),
\end{array}
\] where $\phi_1^0,\phi_2^0\in\mathbb{R}^n$.
\end{thm}
Denoted by $P$ the invariant space of $\mathcal{A}_0$
associated with the eigenvalue zero and $P^*$ the dual space of
$P$, ${\it\Phi}=(\phi_1(\theta),\phi_2(\theta)),-1\leq\theta\leq
0$ and ${\it\Psi}=col(\psi_1(s),\psi_2(s)), 0\leq s\leq1$ the
bases of $P$ and $P^\ast$ correspondingly, we have
$({\it\Psi},{\it\Phi})=I$. Moreover, the following lemma holds.

\begin{lem}{\rm \cite{xh08}}\label{l1}
The bases of $P$ and its dual space $P^*$ have the following
representations:
\[
\begin{array}{lll}
P={\rm span}{\it\Phi},&
{\it\Phi}(\theta)=(\phi_1(\theta),\phi_2(\theta)),&
-1\leq\theta\leq
0,\\
 P^*={\rm span}{\it\Psi},& {\it\Psi}(s)=col(\psi_1(s),\psi_2(s)),&
 0\leq s\leq 1,\end{array}\]
where $\phi_1(\theta)=\phi^0_1\in\mathbb{R}^n\setminus\{0\},
\phi_2(\theta)=\phi^0_2+\phi_1^0\theta, \phi_2^0\in\mathbb{R}^n$
and $\psi_2(s)=\psi^0_2\in\mathbb{R}^{n*}\setminus\{0\},
\psi_1(s)=\psi_1^0-s\psi_2^0, \psi_1^0\in\mathbb{R}^{n*}$, which
satisfy
\begin{equation}\label{e25}
\begin{array}{l}
(1)\ \ (f_1^0+f_2^0)\phi_1^0=0,\\
(2)\ \ (f_1^0+f_2^0)\phi_2^0=(f_2^0+I)\phi_1^0,\\
(3)\ \ \psi_2^0(f_1^0+f_2^0)=0,\\
(4)\ \ \psi_1^0(f_1^0+f_2^0)=\psi_2^0(f_2^0+I),\\
(5)\ \
\psi_1^0\phi_1^0-\frac{1}{2}\psi_2^0f_2^0\phi_1^0+\psi_1^0f_2^0\phi_1^0=1,\\
(6)\ \
\psi_1^0\phi_2^0-\frac{1}{2}\psi_1^0f_2^0\phi_1^0+\psi_1^0f_2^0\phi_2^0+\frac{1}{6}\psi_2^0f_2^0\phi_1^0-\frac{1}{2}\psi_2^0f_2^0\phi_2^0=0,
\end{array}
\end{equation}where  we can determine the unique vector $\phi_1^0,\psi_2^0$
by $(1)$ and $(3)$, respectively, up to some constant factors;
then we can determine $\phi_2^0,\psi_1^0$ by $(2)$ and $(4)$,
respectively. However $(5)$ and $(6)$ are used to determine the
coefficient factors of the vectors $\phi_1^0$ and $\psi_2^0$.
\end{lem}

Based on the characterizations above, we can give another definition for the
quadratic T-B points of (\ref{eyan1}).
\begin{defn}\label{d2}
$(x^0,\lambda^0,\mu^0)$ is called a quadratic T-B point of
(\ref{eyan1}), if it
is a T-B point and\\
\indent${\rm (i)}$ $\psi_2^0f_\lambda^0\neq 0;$\\
\indent${\rm(ii)}$
$d_0=\det\left(\begin{array}{ll}\psi_2^0(A_1+A_2)\phi_1^0&\psi_2^0(B_1+B_2)\phi_1^0\\[2mm]
\begin{array}{l}\psi_1^0(A_1+A_2)\phi_1^0+\\
\psi_2^0(A_1+A_2)\phi_2^0-\\\psi_2^0A_2\phi_1^0
\end{array}
&
\begin{array}{l}\psi_1^0(B_1+B_2)\phi_1^0+\\
\psi_2^0(B_1+B_2)\phi_2^0-\\\psi_2^0B_2\phi_1^0\end{array}\end{array}\right)\neq
0$;\\
\indent${\rm(iii)}$
$\psi_2^0\phi_2^0-\frac{1}{2}\psi_2^0f_2^0\phi_1^0+\psi_2^0f_2^0\phi_2^0\neq
0$,\\
where
\[\begin{array}{ll}A_1=(f_{11}^0+f_{12}^0)\phi_1^0,&
A_2=(f_{21}^0+f_{22}^0)\phi_1^0,\\
B_1=(f_{11}^0+f_{12}^0)\nu_{\lambda\mu}+c_{\lambda\mu}f_{1\lambda}^0+f_{1\mu}^0,&
B_2=(f_{21}^0+f_{22}^0)\nu_{\lambda\mu}+c_{\lambda\mu}f_{2\lambda}^0+f_{2\mu}^0,\end{array}\]
with $c_{\lambda\mu}=-\frac{\psi_2^0f_\mu^0}{\psi_2^0f_\lambda^0}$ and
$\nu_{\lambda\mu}$ satisfying
\[\begin{array}{l}
(f_1^0+f_2^0)\nu_{\lambda\mu}+c_{\lambda\mu}f_\lambda^0+f_{\mu}^0=0,\\
\psi_2^0\nu_{\lambda,\mu}=0.\end{array}
\]
\end{defn}
\begin{thm}\label{teq}
Definition \ref{d2.2} and \ref{d2} are equivalent.
\end{thm}
\proofbegin We only need to show that the corresponding
conditions in Definition \ref{d2.2} and \ref{d2} are equivalent.
We show the conditions (i) in both Definition \ref{d2.2}
and \ref{d2} are equivalent first, that is, $(\psi_2,G_\lambda^0)\neq
0\Leftrightarrow\psi_2^0f_\lambda^0\neq0$. Noting that
\begin{equation}
G_\lambda^0=\left\{\begin{array}{ll}f_\lambda^0,~~~~ &\theta=0,\\
0, &\theta\in[-1,0)\end{array}\right.
\end{equation}
and $\frac{d}{d\lambda}\dot u(\theta)=0$ as $\theta\in[-1,0)$,
we obtain
\begin{equation}
(\psi_2,G_\lambda^0)=\psi_2(0)f_\lambda^0-\int_{-1}^0\int_0^\theta\psi_2(\xi-\theta)[d\eta_{\lambda^0,\mu^0}(\theta)]G_\lambda(u^0,\lambda^0,\mu^0)(\xi)d\xi=\psi_2^0f_\lambda^0.
\end{equation}
Besides,
$(\psi_2,\phi_2)=\psi_2^0\phi_2^0-\frac{1}{2}\psi_2^0f_2^0\phi_1^0+\psi_2^0f_2^0\phi_2^0$,
this confirms that the conditions (iii) in both Definition \ref{d2.2}
and \ref{d2} are equivalent. We only need to show that the
conditions (ii) in both Definition \ref{d2.2} and \ref{d2} are
equivalent next, and it will be a tedious calculation. Firstly, we get
\begin{equation}\label{der1}\begin{array}{ll} (\psi_2,A\phi_1)&
=(\psi_2,G_{uu}^0\phi_1\phi_1)\\&=\psi_2^0(f_{11}^0\phi_1(0)\phi_1
(0)+f_{12}^0\phi_1
(0)\phi_1(-1)\\
&\ \ \ +f_{21}^0\phi_1(-1)\phi_1(0)+f_{22}^0\phi_1(-1)\phi_1(-1))\\
&=\psi_2^0(f_{11}^0\phi_1^0\phi_1^0 +f_{12}^0\phi_1^0\phi_1^0
+f_{21}^0\phi_1^0\phi_1^0
+f_{22}^0\phi_1^0\phi_1^0)\\
&=\psi_2^0(f_{11}^0+f_{12}^0+f_{21}^0+f_{22}^0)\phi_1^0\phi_1^0\\
&=\psi_2^0(A_1+A_2)\phi_1^0.\end{array} \end{equation}

In a similar way, we have
\begin{equation}\label{der2}\begin{array}{ll}
(\psi_2,B\phi_1)&=(\psi_2,G_{uu}^0\nu_{\lambda\mu}\phi_1+c_{\lambda\mu}G_{u\lambda}^0\phi_1+G_{u\mu}^0\phi_1)\\
&=\psi_2^0(f_{11}^0+f_{12}^0+f_{21}^0+f_{22}^0)\nu_{\lambda\mu}\phi_1^0
\\
&\ \ \ +\psi_2^0(f_{1\lambda}^0+f_{2\lambda}^0)c_{\lambda\mu}\phi_{1}^0+\psi_2^0(f_{1\mu}^0+f_{2\mu}^0)\phi_1^0\\
&=\psi_2^0[(f_{11}^0+f_{12}^0+f_{21}^0+f_{22}^0)\nu_{\lambda\mu}\\
&\ \ \ +(f_{1\lambda}^0+f_{2\lambda}^0)c_{\lambda\mu}+(f_{1\mu}^0+f_{2\mu}^0)]\phi_1^0\\
&=\psi_2^0(B_1+B_2)\phi_1^0.
\end{array}
\end{equation}

Since $f_{12}^0=f_{21}^0$, we obtain
\[\begin{array}{ll}
(\psi_2,A\phi_2)&=(\psi_2,G_{uu}^0\phi_1\phi_2)\\&=(\psi_2^0,(G_{uu}^0\phi_1^0)(\phi_2^0+\phi_1^0\theta))\\
&=\psi_2^0[f_{11}^0\phi_1^0(\phi_2^0+\phi_1^0\cdot0)
+f_{12}^0\phi_1^0(\phi_2^0+\phi_1^0\cdot-1)\\&\quad
+f_{21}^0\phi_1^0(\phi_2^0+\phi_1^0\cdot0)
+f_{12}^0\phi_1^0(\phi_2^0+\phi_1^0\cdot-1)]\\&
=\psi_2^0[f_{11}^0\phi_1^0\phi_2^0+f_{12}^0\phi_1^0(\phi_2^0-\phi_1^0)+f_{21}^0\phi_1^0\phi_2^0+f_{22}^0\phi_1^0(\phi_2^0-\phi_1^0)]\\&
=\psi_2^0(f_{11}^0+f_{12}^0+f_{21}^0+f_{22}^0)\phi_1^0\phi_2^0-\psi_2^0(f_{12}^0+f_{22}^0)\phi_1^0\phi_1^0\\&
=\psi_2^0(A_1+A_2)\phi_2^0-\psi_2^0A_2\phi_1^0.
\end{array}
\]
Jointing the above equation with
\[
\begin{array}{ll}(\psi_1,A\phi_1)&=(\psi_1,G_{uu}^0\phi_1\phi_1)\\&
=(\psi_1^0-s\psi_2^0,G_{uu}^0\phi_1^0\phi_1^0)\\&
=\psi_1^0(A_1+A_2)\phi_1^0,
\end{array}
\]
we obtain
\begin{equation}\label{der3}(\psi_2,A\phi_2)+(\psi_1,A\phi_1)=\psi_2^0(A_1+A_2)\phi_2^0-\psi_2^0A_2\phi_1^0+\psi_1^0(A_1+A_2)\phi_1^0.\end{equation}

At last, noting that
\[
\begin{array}{ll}
(\psi_2,B\phi_2)&=(\psi_2,(G_{uu}^0\nu_{\lambda\mu}+c_{\lambda\mu}G_{u\lambda}^0+G_{u\mu}^0)\phi_2)\\
&=(\psi_2^0,G_{uu}^0\nu_{\lambda\mu}(\phi_2^0+\phi_1^0\theta)+c_{\lambda\mu}G_{u\lambda}^0(\phi_2^0+\phi_1^0\theta)+G_{u\mu}^0(\phi_2^0+\phi_1^0\theta))\\
&=\psi_2^0[f_{11}^0\nu_{\lambda\mu}(\phi_2^0+\phi_{1}^0\cdot0)+f_{12}^0\nu_{\lambda\mu}(\phi_2^0+\phi_{1}^0\cdot-1)\\
&\quad+f_{21}^0\nu_{\lambda\mu}(\phi_2^0+\phi_{1}^0\cdot0)+f_{22}^0\nu_{\lambda\mu}(\phi_2^0+\phi_{1}^0\cdot-1)\\
&\quad+c_{\lambda\mu}f_{1\lambda}^0(\phi_2^0+\phi_{1}^0\cdot0)+c_{\lambda\mu}f_{2\lambda}^0(\phi_2^0+\phi_{1}^0\cdot-1)\\
&\quad+f_{1\mu}^0(\phi_2^0+\phi_{1}^0\cdot0)+f_{2\mu}^0(\phi_2^0+\phi_{1}^0\cdot-1)]\\
&=\psi_2^0[f_{11}^0\nu_{\lambda\mu}\phi_2^0+f_{12}^0\nu_{\lambda\mu}(\phi_2^0-\phi_1^0)+f_{21}^0\nu_{\lambda\mu}\phi_2^0+f_{22}^0(\phi_2^0-\phi_1^0)\\
&\quad+c_{\lambda\mu}f_{1\lambda}^0\phi_2^0+c_{\lambda\mu}f_{2\lambda}^0(\phi_2^0-\phi_{1}^0+
f_{1\mu}^0\phi_2^0+f_{2\mu}^0(\phi_2^0-\phi_{1}^0)]\\
&=\psi_2^0(B_1+B_2)\phi_2^0-\psi_2^0B_2\phi_1^0,
\end{array}
\]
and \[\begin{array}{ll}
(\psi_1,B\phi_1)&=(\psi_1,(G_{uu}^0\nu_{\lambda\mu}+c_{\lambda\mu}G_{u\lambda}^0+G_{u\mu}^0)\phi_1)\\
&=\psi_1^0[(f_{11}^0+f_{12}^0+f_{21}^0+f_{22}^0)\nu_{\lambda\mu}\phi_1^0\\
&\quad+c_{\lambda\mu}(f_{1\lambda}^0+f_{2\lambda}^0)\phi_1^0+(f_{1\mu}^0+f_{2\mu}^0)\phi_1^0]\\
&=\psi_1^0(B_1+B_2)\phi_1^0,
\end{array}\]
we have
\begin{equation}\label{der4}\noindent\
(\psi_2,B\phi_2)+(\psi_1,B\phi_1)=\psi_2^0(B_1+B_2)\phi_2^0-\psi_2^0B_2\phi_1^0+\psi_1^0(B_1+B_2)\phi_1^0.
\end{equation}
Collecting Eqs. (\ref{der1}) to (\ref{der4}), we arrive at our assertion.
$\Box$

According to Theorem \ref{t1} and Lemma \ref{l1}, we introduce the
following defining system for the quadratic T-B points of
(\ref{eyan1})

\begin{equation}\label{e1}
H(v)=\left(\begin{array}{l}
f(x,x,\lambda,\mu),\\[3mm]
(f_1(x,x,\lambda,\mu)+f_2(x,x,\lambda,\mu))\varphi_1\\[3mm]
(f_1(x,x,\lambda,\mu)+f_2(x,x,\lambda,\mu))\varphi_2-(f_2(x,x,\lambda,\mu)+I)\varphi_1\\[3mm]
l_1\varphi_1-\frac{1}{2}l_2f_2(x,x,\lambda,\mu)\varphi_1+l_1f_2(x,x,\lambda,\mu)\varphi_1-1\\[3mm]
l_1\varphi_2-\frac{1}{2}l_1f_2(x,x,\lambda,\mu)\varphi_1+l_1f_2(x,x,\lambda,\mu)\varphi_2+\\\hspace{20mm}
\frac{1}{6}l_2f_2(x,x,\lambda,\mu)\varphi_1-
\frac{1}{2}l_2f_2(x,x,\lambda,\mu)\varphi_2
\end{array}\right)=0,
\end{equation}
where $v=(x,\varphi_1,\varphi_2,\lambda,\mu)^T\in
V=\mathbb{R}^n\times\mathbb{R}^n\times\mathbb{R}^n\times\mathbb{R}\times\mathbb{R}$,
$H:V\rightarrow
V$. From the discussions above, we know the system is
equivalent to (\ref{ea2}).

\begin{thm}\label{t2}
Assume $(x^0,\lambda^0,\mu^0)$ is a quadratic T-B point of
(\ref{eyan1}), then the defining system (\ref{e1}) is regular at its zero
$v^0=(x^0,\phi_1^0,\phi_2^0,\lambda^0,\mu^0)^T$.
\end{thm}
\proofbegin
 Obviously, $H_v^0$ reads
{\footnotesize
$$
\left(\begin{array}{lllll} f_1^0+f_2^0&0&0&f_\lambda^0&f_\mu^0\\[6mm]
(f_{11}^0+f_{12}^0+f_{21}^0+f_{22}^0)\phi_1^0&f_1^0+f_2^0&0&(f_{1\lambda}^0+f_{2\lambda}^0)\phi_1^0&(f_{1\mu}^0+f_{2\mu}^0)\phi_1^0\\[6mm]
\begin{array}{l}(f_{11}^0+f_{12}^0+f_{21}^0+f_{22}^0)\phi_2^0-\\[1mm](f_{21}^0+f_{22}^0)\phi_1^0\end{array}&-(f_2^0+I)&f_1^0+f_2^0&\begin{array}{l}(f_{1\lambda}^0+f_{2\lambda}^0)\phi_2^0-\\f_{2\lambda}^0\phi_1^0\end{array}&\begin{array}{l}(f_{1\mu}^0+f_{2\mu}^0)\phi_2^0-\\f_{2\mu}^0\phi_1^0\end{array}\\[6mm]
\begin{array}{l}-\frac{1}{2}l_2(f_{12}^0+f_{22}^0)\phi_1^0+\\l_1(f_{12}^0+f_{22}^0)\phi_1^0\end{array}&l_1-\frac{1}{2}l_2f_2^0+l_1f_2^0&0&\begin{array}{l}-\frac{1}{2}l_2f_{1\lambda}^0\phi_1^0+\\l_1f_{2\lambda}^0\phi_1^0\end{array}&\begin{array}{l}-\frac{1}{2}l_2f_{2\mu}^0\phi_1^0+\\l_1f_{2\mu}^0\phi_1^0\end{array}\\[6mm]
\begin{array}{l}-\frac{1}{2}l_1(f_{12}^0+f_{22}^0)\phi_1^0+\\l_1(f_{12}^0+f_{22}^0)\phi_2^0+\\
\frac{1}{6}l_2(f_{12}^0+f_{22}^0)\phi_1^0-\\\frac{1}{2}l_2(f_{12}^0+f_{22}^0)\phi_2^0
\end{array}
&-\frac{1}{2}l_1f_2^0+\frac{1}{6}l_2f_{2}^0&\begin{array}{l}l_1+l_1f_2^0-\\\frac{1}{2}l_2f_2^0\end{array}
&\begin{array}{l}-\frac{1}{2}l_1f_{2\lambda}^0\phi_1^0+\\
l_1f_{2\lambda}^0\phi_2^0+\\
\frac{1}{6}l_2f_{2\lambda}^0-\\
\frac{1}{2}l_2f_{2\lambda}^0\phi_2^0
\end{array}
&\begin{array}{l}-\frac{1}{2}l_1f_{2\mu}^0\phi_1^0+\\
l_1f_{2\mu}^0\phi_2^0+\\[1mm]
\frac{1}{6}l_2f_{2\mu}^0-\\[1mm]
\frac{1}{2}l_2f_{2\mu}^0\phi_2^0
\end{array}
\end{array}\right).
$$
} We first prove the map $H_v^0$ is injective. By expanding
 $H_v^0\vartheta=0$ with
 \ $\vartheta=(\vartheta_1,\vartheta_2,\vartheta_3,c_1,c_2)^T\in V$, we obtain
\begin{equation}\label{bs1}
(f_1^0+f_2^0)\vartheta_1+c_1f_\lambda^0+c_2f_\mu^0=0,
\end{equation}
\begin{equation}\label{bs2}
(f_{11}^0+f_{12}^0+f_{21}^0+f_{22}^0)\phi_1^0\vartheta_1+(f_1^0+f_2^0)\vartheta_2+(f_{1\lambda}^0+f_{2\lambda}^0)\phi_1^0c_1+(f_{1\mu}^0+f_{2\mu}^0)\phi_1^0c_2=0,
\end{equation}
\begin{equation}\label{bs3}
\begin{array}{l}\Big((f_{11}^0+f_{12}^0+f_{21}^0+f_{22}^0)\phi_2^0-(f_{21}^0+f_{22}^0)\phi_1^0\Big)\vartheta_1-(f_2^0+I)\vartheta_2+(f_1^0+f_2^0)\vartheta_3\\[2mm]
\qquad\qquad\quad+\Big((f_{1\lambda}^0+f_{2\lambda}^0)\phi_2^0-f_{2\lambda}^0\phi_1^0\Big)c_1+
\Big((f_{1\mu}^0+f_{2\mu}^0)\phi_2^0-f_{2\mu}^0\phi_1^0\Big)c_2=0,\end{array}
\end{equation}
\begin{equation}\label{bs4}
\begin{array}{l}\Big(-\frac{1}{2}l_2(f_{12}^0+f_{22}^0)\phi_1^0+l_1(f_{12}^0+f_{22}^0)\phi_1^0\Big)\vartheta_1+(l_1-\frac{1}{2}l_2f_2^0+l_1f_2^0)\vartheta_2+\\[2mm]
\qquad\qquad\quad+(-\frac{1}{2}l_2f_{1\lambda}^0\phi_1^0+l_1f_{2\lambda}^0\phi_1^0)c_1+(-\frac{1}{2}l_2f_{2\mu}^0\phi_1^0+l_1f_{2\mu}^0\phi_1^0)c_2=0
\end{array}
\end{equation}
and
\begin{equation}\label{bs5}
\begin{array}{l}\Big(-\frac{1}{2}l_1(f_{12}^0+f_{22}^0)\phi_1^0+l_1(f_{12}^0+f_{22}^0)\phi_2^0\\
\qquad\qquad\quad+\frac{1}{6}l_2(f_{12}^0+f_{22}^0)\phi_1^0-\frac{1}{2}l_2(f_{12}^0+f_{22}^0)\phi_2^0
\Big)\vartheta_1\\[2mm]
\qquad\qquad\quad+(-\frac{1}{2}l_1f_2^0+\frac{1}{6}l_2f_{2}^0)\vartheta_2+(l_1+l_1f_2^0-\frac{1}{2}l_2f_2^0)\vartheta_3\\[2mm]
\qquad\qquad\quad+(-\frac{1}{2}l_1f_{2\lambda}^0\phi_1^0+l_1f_{2\lambda}^0\phi_2^0+
\frac{1}{6}l_2f_{2\lambda}^0-\frac{1}{2}l_2f_{2\lambda}^0\phi_2^0)c_1\\[2mm]
\qquad\qquad\quad+(-\frac{1}{2}l_1f_{2\mu}^0\phi_1^0+l_1f_{2\mu}^0\phi_2^0+
\frac{1}{6}l_2f_{2\mu}^0-\frac{1}{2}l_2f_{2\mu}^0\phi_2^0)c_2=0.
\end{array}
\end{equation}
To multiply (\ref{bs1}) by $\psi_2^0$ from left, we have by Lemma
\ref{l1}
\[
\psi_2^0f_\lambda^0c_1+\psi_2^0f_\mu^0c_2=0,
\]
which yields
$c_1=-\frac{\psi_2^0f_\mu^0}{\psi_2^0f_\lambda^0}c_2$. Noting that
$c_{\lambda\mu}=-\frac{\psi_2^0f_\mu^0}{\psi_2^0f_\lambda^0}$, we
obtain $c_1=c_{\lambda\mu}c_2$, therefore,
$\vartheta_1=c\phi_1^0+c_2\nu_{\lambda\mu}$, where $c$ is a
constant to be determined. By substituting them into Eqs.
(\ref{bs2}) and (\ref{bs3}), it yields
\begin{equation}\label{bs11}
c(A_1+A_2)\phi_1^0+c_2(B_1+B_2)\phi_1^0+(f_1^0+f_2^0)\vartheta_2=0,
\end{equation}
and
\begin{equation}\label{bs21}
c(A_1+A_2)\phi_2^0+c_2(B_1+B_2)\phi_2^0-cA_2\phi_1^0-c_2B_2\phi_1^0-(f_2^0+I)\vartheta_2+(f_1^0+f_2^0)\vartheta_3=0.
\end{equation}
To multiply (\ref{bs11}) by $\psi_2^0$ from left, we obtain by
Lemma \ref{l1} again
\begin{equation}\label{tri}
c\psi_2^0(A_1+A_2)\phi_1^0+c_2\psi_2^0(B_1+B_2)\phi_1^0=0.
\end{equation}
To multiply (\ref{bs11}) by $\psi_1^0$ from left,  we have
\begin{equation}\label{lin1}
c\psi_1^0(A_1+A_2)\phi_1^0+c_2\psi_1^0(B_1+B_2)\phi_1^0+\psi_1^0(f_1^0+f_2^0)\vartheta_2=0.
\end{equation}
Multiplying  (\ref{bs21}) by $\psi_2^0$ from left yields
\begin{equation}\label{lin2}
\begin{array}{l}
c\psi_2^0(A_1+A_2)\phi_2^0+c_2\psi_2^0(B_1+B_2)\phi_2^0-c\psi_2^0A_2\phi_1^0\\\quad-c_2\psi_2^0B_2\phi_1^0-\psi_2^0(f_2^0+I)\vartheta_2+\psi_2^0(f_1^0+f_2^0)\vartheta_3=0.
\end{array}
\end{equation}
Adding (\ref{lin1}) and (\ref{lin2}) together, and utilizing
(\ref{e25}) we obtain
\begin{equation}\label{ead}
\begin{array}{l}
c(\psi_1^0(A_1+A_2)\phi_1^0+\psi_2^0(A_1+A_2)\phi_2^0-\psi_2^0A_2\phi_1^0)\\\quad+c_2(\psi_1^0(B_1+B_2)\phi_1^0+\psi_2^0(B_1+B_2)\phi_2^0-\psi_2^0B_2\phi_1^0)=0.
\end{array}
\end{equation}

From (\ref{tri}) and (\ref{ead}), noting that $d_0\neq 0$, we
obtain $c=c_2=0$, therefore $c_1=0$, $\vartheta_1=0$.
Consequently, (\ref{bs11}) reads $(f_1^0+f_2^0)\vartheta_2=0$,
therefore, $\vartheta_2\in\mathcal{N}(f_1^0+f_2^0)$. If
$\vartheta_2\neq0$, we have
\[
l_1\vartheta_2-\frac{1}{2}l_2f_2^0\vartheta_2+l_1f_2^0\vartheta_2=(l,\vartheta_2)\neq0,
\]
which contradicts with (\ref{bs4}), hence $\vartheta_2=0$. Hereby,
(\ref{bs3}) reads $ (f_1^0+f_2^0)\vartheta_3=0$, together with
(\ref{bs5}) it implies $\vartheta_3=0$. Altogether $\vartheta=0$.

In a similar way, we can prove $H_v^0$ is surjective. This concludes the theorem.
$\Box$

Noting that the system (\ref{e1}) is of finite dimension, according to Theorem \ref{t2}, it can be solved by many iteration procedures, especially the
standard Newton iteration procedure, that is
\begin{equation}\label{Newton}
v_{k+1}=v_{k}-[J_H(v_k)]^{-1}H(v_k) , \ \mbox{for}\  k\geqslant 0
\end{equation}
where $J_H(\cdot)$ is the Jacobian matrix.

\section{Numerical example}
Consider the following predator prey system with delay \cite{ZhL05},
\begin{eqnarray}\label{s1}
\left\{ \begin{array}{ll}
\dot x_1(t)=rx_1(t)(1-\frac{x_1(t)}{K})-\frac{x_1(t-\tau)x_2(t)}{a+x_1^2(t-\tau)}, \\
\dot x_2(t)=x_2(t)(\frac{\mu x_1(t-\tau)}{a+x_1^2(t-\tau)}-D)
\end{array}
\right.
\end{eqnarray}
where $r, K, a, \mu, D$ and $\tau $ are positive constants. Eq.
(\ref{s1}) has a T-B singularity if the parameters $\mu, a, D$ and
$K$ satisfy $\mu^2-4aD^2=0$ and $\mu=KD$ \cite{ZhL05}. For
determining the T-B point of (\ref{s1}), we choose $K$ and $D$ as
parameters and fix the others. Specially, taking $a=1$, $\mu =1$,
$\tau =1$, and $r=1$ in (\ref{s1}), we have
\begin{equation}\label{add}
\dot x(t)=f(x(t),x(t-1),D,K)
\end{equation}
with $x(t)=(x_1(t),x_2(t))^T$ and
\[f(x(t),x(t-1),D,K)=\left(\begin{array}{c}
x_1(t)(1-\frac{x_1(t)}{K})-\frac{x_1(t-1)x_2(t)}{1+x_1^2(t-1)}\\
x_2(t)(\frac{ x_1(t-1)}{1+x_1^2(t-1)}-D)
\end{array}\right).\]
For any given $x\in \mathbb{R}^2$, the derivatives of $f$ with
respect to the first and the second variables at $x$ are
respectively given by

\[
f_1(x,x,D,K)=\left(\begin{array}{cc}
1-\frac{2x_1}{K} &-\frac{x_1}{1+x_1^2}\\[4mm]
0&\frac{x_1}{1+x_1^2}-D\end{array}\right)
\]
and
\[
f_2(x,x,D,K)=\left(\begin{array}{cc}
-x_2\frac{(1-x_1^2)}{(1+x_1^2)^2} & 0\\[4mm]
x_2\frac{(1-x_1^2)}{(1+x_1^2)^2}&0\end{array}\right).\]
Hence the defining system for the quadratic T-B point of
(\ref{s1}) reads
\begin{equation}\label{ex1}
H(v)=\left(\begin{array}{l}
f(x,x,D,K)\\[3mm]
(f_1(x,x,D,K)+f_2(x,x,D,K))\varphi_1\\[3mm]
(f_1(x,x,D,K)+f_2(x,x,D,K))\varphi_2-(f_2(x,x,D,K)+I)\varphi_1\\[3mm]
l_1\varphi_1-\frac{1}{2}l_2f_2(x,x,D,K)\varphi_1+l_1f_2(x,x,D,K)\varphi_1-1\\[3mm]
l_1\varphi_2-\frac{1}{2}l_1f_2(x,x,D,K)\varphi_1+l_1f_2(x,x,D,K)\varphi_2+\\\hspace{20mm}
\frac{1}{6}l_2f_2(x,x,D,K)\varphi_1-
\frac{1}{2}l_2f_2(x,x,D,K)\varphi_2
\end{array}\right)=0,
\end{equation}
where $v=(x,\varphi_1,\varphi_2,D,K)^T\in
V=\mathbb{R}^2\times\mathbb{R}^2\times\mathbb{R}^2\times\mathbb{R}\times\mathbb{R}$,
    $H:V\rightarrow V$.

    Choosing $l_1=(1,0)$ and $l_2=(1,0)$ and applying Newton method
(\ref{Newton}) to Eq. (\ref{ex1}), we obtain the following results
shown in Table 1 with Matlab. In each iteration with respect to different initial values, the solution converges to $(1,1,1,0,0,-2,0.5,2)$ with the remainder less than $10^{-21}$.
\begin{center}
\begin{tabular}{|l|c|}
\hline
Initial value $v_0$ &Step \\
\hline
(1.1,   1.1,  1,   0, 3, 0, 0.4,  1) & 5\\
 \hline
 (1.2, 1.2, 1.2,  1, 1,  0,  0.5, 0.5)& 7\\
 \hline
 (1.5, 1.5, 1.5,  1.5, 1.5,  1.5,  0.6,  1.6)& 7\\
 \hline
 (3, 1.5, 1.2, 0.5, 1.8, -1.8,  0.45,  1.9) & 6\\
 \hline
\end{tabular}\\
\vspace*{2mm} Table 1:\ {\small Number of iterations required by Newton methods for (\ref{ex1}) with respect to different initial values.}
\end{center}

%
%

In fact, according to \cite{ZhL05} we can obtain exactly  the T-B point of (\ref{s1}). For $a=\mu=\tau=r=1$, the
T-B point is $(x_1,x_2)=(1,1)$ with parameter values
$D=\frac{1}{2}$ and $K=2$.

\begin{remark}From the numerical result we see that the Newton's method converges
very rapidly once the initial value is close to the true solution. However, choosing the initial point $v_0$ to
ensure a convergent Newton iteration for solving
(\ref{e1}) is a rather complicated problem because the function
$H(v)$ is highly nonlinear and high dimensional. There is
no good method as so far even for ODEs or for the first order singular points without the aid of continuation techniques.
The guess of the initial value $v_0$ for
starting the iteration procedure can be found in the following
way. First we fix one parameter, for instance $\mu=\mu^H$. According to the method by Luzyanina and Roose \cite{Dr96}, we can find a Hopf
bifurcation point $x^H$ of (\ref{eyan1}) on its Hopf bifurcation
curve at $\lambda=\lambda^H$, as well the approximation of the eigenvector
$\phi_{Re}^H+i\phi_{Im}^H$ associated with eigenvalue $i\omega^H$.
Next, we start from $(x^H,\lambda^H,\mu^H,\omega^H,\phi_{Re}^H,\phi_{Im}^H)$ to trace the Hopf bifurcation curves by the continuation technique. When the singularity of the nonlinear system used for
continuation changed at some point, for
instance $(x_0,\lambda_0,\mu_0,\omega_0,\phi_{Re}^0,\phi_{Im}^0)$,
then we can use $(x_0,\phi_{Re}^0,\bar \phi,\lambda_0,\mu_0)$ as initial value to
start the Newton iteration procedure, where $\bar \phi$ still
should be guessed, but it is not a hard work any more.
\end{remark}











\end{document}